\documentclass[a4paper,12 pt,titlepage]{article}
\usepackage{amsmath, amsthm, amssymb
}
\usepackage{multirow}
\usepackage{booktabs}
\usepackage{mathrsfs}
\usepackage{verbatim}
\usepackage[pdftex]{graphicx}
\newtheorem{theorem}{Theorem}
\newtheorem{lemma}{Lemma}
\newtheorem{prop}{Property}

\theoremstyle{definition}
\newtheorem{definition}{Definition}
\usepackage{geometry}
\usepackage[pdftex]{hyperref}
\hypersetup{
colorlinks=true, linkcolor=blue, bookmarks=true}
\title{Modeling Air Traffic Flow Management}

\author{Hrishikesh Ganu\\PGP Student\\PGP-10-80\\
\texttt{hrishikeshv10@iimk.ac.in}\\ \\ \\ \emph{Under the guidance of} \\ \\ 
Dr. Pritibhushan Sinha\\
Assistant Professor\\Quantitative Methods and Operations Management Area \\Indian Institute of Management, Kozhikode\\ P.O. Kozhikode-673570\\
 \texttt{pritibhushan.sinha@iimk.ac.in}}

\begin{document}
\maketitle

\setcounter{page}{1}
\pagestyle{myheadings}
\markright{\textsc{Modeling Air Traffic Flow Management}}

\pagenumbering{roman}
\tableofcontents 
\newpage
\listoffigures 
\listoftables 
\newpage
\section*{Acknowledgments}
I am grateful to my guide, Dr. Pritibhushan Sinha, for his constant support and encouragement. It is due to him, that I could work not only on the theoretical but also on the practical parts of the problem. The procedure for reducing the problem size, mentioned in the ``Directions for future work,'' is a result of his valuable suggestions. This course was indeed a refreshing change from the otherwise descriptive courses at IIM-K.  I have benefited a lot from discussions with Dr. Saji Gopinath, regarding puzzles and other applications of Integer Programming, which would be reflected implicitly through this work.  

I would like to thank Ms. T.Sunitha,  Assistant Librarian at the Library and Information Centre of IIM-Kozhikode for help in procuring  reference materials. I was previously a graduate student in Aerospace Engineering at IISc Bangalore and the initial idea for this project was to some extent, a result of discussions with my guide, Professor B.N. Raghunandan. The GNU Mathprog language that we have used here is available freely only because of the Open Source movement and thanks are due to  Andrew Makhorin and others who develop and maintain it.

\pagenumbering{arabic}
\setcounter{page}{1}

\newpage
\section{The problem}
Air Traffic Flow Management is the regulation of air traffic in order to avoid exceeding 
airport or flight sector capacity in handling traffic, and to ensure that available capacity is used efficiently. The number of 
flights taking off or landing from a certain airport or the number of planes traveling in a particular sector are functions
of several variables including:
\begin{enumerate}
\item \textbf{The number of runways available}
\item\textbf{ATC capacity}\\The Air Traffic Controllers and the Control centre perform various tasks like:
 Terminal Approach,
Final Approach,
Ground (arrivals),
Gate operations (arrivals),
Airport Surface,
Gate operations (departures),
Ground (departures) and finally the
Take-Off / Transition to Terminal and Center. Thus the number and expertise of ATCs available determines the airport capacity to a certain extent.
\item\textbf{Restrictions as to which aircraft can follow an aircraft of a given class}\\
This is because an aircraft periodically sheds vortices from the wing tips. These vortices can interact with the boundary layer
around the wing of an aircraft following in the first aircraft's wake and destroy the lift. The effect is more pronounced
 when the following aircraft is of a smaller size class then the leading one. This means that a delay must be imposed such that
$F_{t}\geq L_{t}+\delta_{F,L}$ where $F_{t}$ and $L_{t}$ denote the times  at which the follower and the leader take-off or 
land and $\delta_{F, L}$ denotes the delay.
\item\textbf{Airspace restrictions}\\ Aircraft have to follow certain ``corridors''  when traveling in air. This restricts the
number of aircraft that can be airborne at any given instant.
Finally, most of  the restrictions mentioned above, are functions of time too. This makes the problem dynamic in nature.
\end{enumerate}
\section{Existing models}

Airport and airspace capacity is the major
cause of congestion. The difficulty in dealing with
this parameter is its uncertainty heavily influenced
by weather conditions among other factors.
 Short-term solutions (12 hours time horizon) for air
traffic flow management include ground-holding
policies. These policies are motivated by the fundamental
fact that airborne delays are much
costlier than ground delays, since the former include
fuel, maintenance, depreciation and safety
costs. Therefore, the aim of ground-holding policies
is to translate anticipated airborne delays to
the ground.
In addition to ground holding, distributing the traffic efficiently across the airspace helps in reducing congestion
in the network.

Several models have been proposed in the literature
for problem solving under different real-life
hypotheses. The Single-Airport Ground 
Problem (SAGHP) only considers a single airport
and the goal is to produce ground-holding schedules. 
 The Multi-Airport Ground-Holding Problem
(MAGHP), Andreatta et.al.\cite{Andreatta} considers a network of airports such
that the ground-holding policies for one of them
have impact on the other airport schedules, as can be seen in Vranas et.al. \cite{Vranas}
. The Air Traffic Flow Management Problem
(TFMP) is similar to MAGHP but also considers the
airspace network besides the airport capacity,
Bertsimas et.al. \cite{Bert}. The main assumption on the problem
made by most of the approaches described in the
open literature is the deterministic character of the
data.
\section{Motivation and Objectives}
The present work focuses on the seminal paper, due to Bertsimas et.al \cite{Bert} in which the polyhedral structure of the TFMP is analyzed. The authors of this paper claim that by  using a  transformation of the decision variables, a formulation in which some inequalities  define the facets of the convex hull of the TFMP was obtained.

This is a significant improvement in the formulation, since Polyhedral theory indicates that under the conditions mentioned above,
 the LP relaxation of the  TFMP, might be integral. This means that we may not need techniques like branch and bound, at least in some of the cases.
In this work an attempt has been made to understand and verify the claims made in the paper \cite{Bert}. To be specific:
\begin{enumerate}
\item The physical problem was studied in terms of the several  
formulations\\ (SAGHP, MAGHP, TFMP) available.
\item We tried to understand the main results of of Polyhedral Combinatorics, principally through the textbook by Wolsey and Nemhauser~\cite{Wolsey}.
\item The paper due to Bertsimas et.al \cite{Bert} was examined in the light of this theory to understand how the formulation helps
to achieve integral solutions.
\end{enumerate}
\section{Method}
\subsection{The Traffic Flow Management Problem (TFMP)}
This subsection has been reproduced verbatim from Bertsimas et.al~\cite{Bert}.
 Consider a set of flights, $\mathscr{ F} = \{1, ... ,F\}$, a set of airports, $\mathscr{K} = \{1, ... ,K\}$, a set of time periods,
$\mathscr{T} = \{1, ... ,T\}$, and a set of pairs of flights that are continued, $\mathscr{C} = \{(\acute{f}, f) : \acute{f}\textrm{ is continued
by flight} f\}$. We shall refer to any particular time period $t$ as the ``time $t$''. The problem input
data are given as follows:\\
$N_{f}$ = number of sectors in flight $f's$ path\\\\
$P(f, i)$ = the $i^{th}$ sector in flight $f's$ path\\\\
$P_{f} = \{P(f,i): 1 \leq i \leq N_{f}\}$\\\\
$D_{k}(t) $= departure capacity of airport $k$ at time $t$\\\\
$A_{k}(t) $= arrival capacity of airport $k$ at time $t$\\\\
$S_{j}(t)$ = capacity of sector $j$ at time $t$\\\\
$d_{f}= $ scheduled departure time of flight $f$\\\\
$r_{f}=$  scheduled arrival time of flight $f$\\\\
$s_{f}=$  turnaround time of an airplane after flight $f$\\\\
$c_{g}^{f}=$  cost of holding flight $f$ on the ground for one unit of time\\\\
$c_{a}^{f} =$ cost of holding flight $f$ in the air for one unit of time\\\\
$l_{fj} =$ number of time units that flight $f$ must spend in sector $j$\\\\
$T_{f}^{j} =$ set of feasible times for flight $f$ to be in sector $j$\\\\
Note that by ``flight,'' we mean a ``flight leg'' between two airports. Also, flights referred
to as ``continued'' are those flights whose aircraft is scheduled to perform a later flight within
some time interval of its scheduled arrival.\\
\textbf{Objective}: The objective in the TFMP is to decide how much each flight is going to be held
on the ground and in the air in order to minimize the total delay cost.
We model the problem as follows.\\
\textbf{Decision variables:}
\begin{displaymath}
w^{j}_{ft}= \left\{ \begin{array}{ll}
1 & \textrm{if flight $f$ arrives at sector $j$ by time $t$}\\
0 & \textrm{otherwise}
\end{array} \right.
\end{displaymath}
 
Note that the $w^{j}_{ft}$ are defined as being 1 if flight $f$ arrives at sector $j$ by time $t$. This definition
using \textit{by} and not \textit{at} is critical to the understanding of the formulation. Also recall that we have
also defined for each flight a list $P_{f}$ of sectors which includes the departure and arrival airports,
so that the variable $w^{j}_{ft}$ will only be defined for those sectors $j$ in the list $P_{f}$. Moreover, we
have defined $T^{j}_{f}$ as the set of feasible times for flight $f$ to be in sector $j$, so that the variable $w^{j}_{ft}$ will only be defined for those times within  $T^{j}_{f}$. Thus, in the formulation whenever the variable
$w^{j}_{ft}$ is used, it is assumed that this is a feasible $(f, j, t)$ combination. Furthermore, one variable
per flight-sector pair can be eliminated from the formulation by setting $w^{j}_{f\bar{T_{f}^{j}}}$ = 1 where $\bar{T_{f}^{j}}$ is the last time period in the set $T_{f}^{j}$ . Since flight $f$ has to arrive at sector $j$ by the last possible
time in its time window, we can simply set it equal to one as a parameter before solving the
problem.

Having defined the variables  we can express several quantities of interest as linear
functions of these variables as follows:
\begin{enumerate}
\item Noticing that the first sector for every flight represents the departing airport, then the
total number of time units that flight $f$ is held on the ground is the actual departure
time minus the scheduled departure time, i.e.,
\begin{displaymath}
g_{f}=\sum_{t \in T^{k}_{f}, k=P(f,1)}t(w^{k}_{ft}-w^{k}_{f,t-1}) - d_{f}\end{displaymath}.
\item Noticing that the last sector for every flight represents the destination airport, the total
number of time units that flight $f$ is held in the air can be expressed as the actual arrival
time minus the scheduled arrival time minus the amount of time that the flight has been
held on the ground, i.e.,
$$
a_{f}=\sum_{t \in T^{k}_{f}, k=P(f,N_{f})}t(w^{k}_{ft}-w^{k}_{f,t-1}) -r_{f} -g_{f}
$$
\end{enumerate}
The objective of the formulation is to minimize total delay cost, and the TFMP is hence:
$$
\mathbf{TFMP}:\\
\textrm{Min} \quad z=\sum_{f\in \mathscr{F}} c^{g}_{f} g_{f}+ c^{a}_{f}a_{f}
$$

\begin{equation}
\textrm{subject to} \sum_{f:P(f,1)=k)}(w^{k}_{ft}-w^{k}_{f,t-1}) \leq D_{k}(t),\forall
 k \in \mathscr{K}, t \in \mathscr{T}
\label{depcap}
\end{equation}
\begin{equation}
 \sum_{f:P(f,N_{f})=k)}(w^{k}_{ft}-w^{k}_{f,t-1}) \leq A_{k}(t),\forall
 k \in \mathscr{K}, t \in \mathscr{T}
\label{arrvcap}
\end{equation}
\begin{equation}
 \sum_{f:P(f,i)=j,P(f,i+1)=\acute{j}}(w^{j}_{ft}-w^{\acute{j}}_{f,t}) \leq S_{j}(t),\forall
 j \in \mathscr{J}, t \in \mathscr{T}
\label{sectorcap}
\end{equation}

\begin{equation}
w^{\acute{j}}_{f,t+l_{fj}}-w^{j}_{f,t} \leq 0,\left\{
\begin{array}{ll}
\forall f \in \mathscr{F}, t \in T_{f}^{j},j=P(f,i),\\
\acute{j}=P(f, i+1), i < N_{f}
\end{array}\right.
\label{transit}
\end{equation}

\begin{equation}
w^{k}_{f,t}-w^{k}_{\acute{f},t-s_{\acute{f}}} \leq 0,\left\{
\begin{array}{ll}
\forall (f, \acute{f}) \in \mathscr{C}, t \in T_{f}^{k},\\
P(f,1)=k=P(\acute{f}, N_{\acute{f}})
\end{array}\right.
\label{turn}
\end{equation}

\begin{equation}
w^{j}_{f,t}-w^{j}_{f,t-1} \geq 0,\forall f \in \mathscr{F},
 j \in P_{f}, t \in \mathscr{T_{f}^{j}}
\label{by}
\end{equation}

\begin{equation}
w^{j}_{f,t} \in \{0,1\}, \forall f \in \mathscr{F},
 j \in P_{f}, t \in \mathscr{T_{f}^{j}}
\label{binary}
\end{equation}

The first three constraints take into account the capacities of various aspects of the system.
The  constraint (\ref{depcap}) ensures that the number of flights which may take off from airport $k$ at
time t, will not exceed the departure capacity of airport $k$ at time $t$. Likewise, the constraint (\ref{arrvcap}) ensures that the number of flights which may arrive at airport $k$ at time $t$, will not
exceed the arrival capacity of airport $k$ at time $t$. In each case, the difference will be equal to
one only when the first term is one and the second term is zero. Thus, the differences capture
the time at which a flight uses a given airport. The constraint (\ref{sectorcap}) ensures that the sum of
all flights which may feasibly be in sector $j$ at time $t$ will not exceed the capacity of sector $j$
at time $t$. This difference gives the flights which are in sector $j$ at time $t$, since the first term
will be 1 if flight $f$ has arrived in sector $j$ by time $t$ and the second term will be 1 if flight $f$
has arrived at the next sector by time $t$. So the only flights which will contribute a value 1 to
this sum are the flights that have arrived at $j$ and not yet departed by time $t$.
Constraints (\ref{transit}) represent connectivity between sectors. They stipulate that if a flight
arrives at sector $\acute{j}$ by time $t + l_{fj}$, then it must have arrived at sector $j$ by time $t$ where $j$
and $\acute{j}$ are contiguous sectors in flight $f's$ path. In other words, a flight cannot enter the next
sector on its path until it has spent $l_{fj} $time units (the minimum possible) traveling through
sector $j$, the current sector in its path.
Constraints (\ref{turn}) represent connectivity between airports. They handle the cases in which a
flight is continued, i.e., the flight's aircraft is scheduled to perform a later flight within some
time interval. We will call the first flight $\acute{f}$ and the following flight $f$. Constraints (\ref{turn}) state
that if flight $f$ departs from airport $k$ by time $t$, then flight $\acute{f}$ must have arrived at airport $k$
by time $t - s_{\acute{f}}$. The turnaround time, $s_{\acute{f}}$, takes into account the time that is needed to clean,
refuel, unload and load and further prepare the aircraft for the next flight. In other words,
flight $f$ cannot depart from airport $k$, until flight $\acute{f}$ has arrived and spent at least $s_{\acute{f}}$ time
units at airport $k$.
Constraints (\ref{by}) represent connectivity in time. Thus, if a flight has arrived by time $t$, then
$w_{ft}$, has to have a value of 1 for all later time periods, $\acute{t} \geq t$.
\subsubsection{Complexity of the TFMP}
Bertsimas et.al.,~\cite{Bert}, proved that the TFMP with all capacities equal to 1 is NP-hard, by showing equivalence with the Job-shop scheduling problem.
\subsubsection{Computational results due to Bertsimas et.al.~\cite{Bert}}
Bertsimas et.al, solved many instances of the MAGHP and the TFMP using a CPLEX 2.1 solver. They used as many as 1000 flights, varying the number of connected flights from 1/5 to 4/5 of the total flights being flown. The problem was solved over a 24-hour period, with intervals of 15 min. Even at the infeasibility border (which represents the data set such that even a minor change in parameters will make the constraints inconsistent), only 4\% of the solutions were found to be non-integral, that too for the MAGHP. All instances of the TFMP had integral solutions.

\section{Polyhedral combinatorics to identify Facet defining constraints }
The use of the variables denoting whether a flight has arrived \textit{by} a certain time rather than \textit{at} a certain time has advantages, achieved in terms of a ``stronger'' formulation. Key results from Polyhedral theory, relevant to the present work are included here and are used in sketching a constructive proof of the high incidence of integral solutions even for the LP relaxation of TFMP.
\subsection{Basic Concepts from Polyhedral Theory }

The  definitions and theorems in this section, closely follow the development given in Wolsey and Nemhauser \cite{Wolsey}.
\begin{definition}
A set of $k$ vectors,$\in \Re^{n}$  is defined to be affinely independent if \begin{displaymath}
\sum_{i=1}^k \lambda_{i} x^{i}=0 ,  \sum_{i=1}^{k}\lambda_{i}=0 \Rightarrow  \lambda_{i}=0 ~  \forall i=1,....,k \end{displaymath}
\end{definition}
\begin{definition}[\textbf{Polyhedron}]
A \textit{polyhedron} $P\subseteq\Re^{n}$ is a set of points that satisfy a finite number of linear inequalities; that is, $P=\{x\in \Re^{n}: Ax\leq b \}$, where $(A,b)$ is an $m\times(n+1)$ matrix.\end{definition}
\begin{definition}\label{dimension}
The \textit{dimension} of a Polyhedron $P$ is one less than the cardinality of the maximal set of affinely independent points in $P$.
\end{definition}
\begin{lemma}[\textbf{Dimension of $P$}]
If $P\subseteq \Re^{n}$, then $dim(P)+rank(A^{=},b^{=})=n$ where $(A^{=},b^{=})$ is the equality set consisting of
 the rows corresponding to $a^{i}x=b_{i} 
\quad \forall x\in P$.
\end{lemma}
\begin{definition}[\textbf{Valid Inequality}]
An inequality $\pi x \leq \pi_{0}$ is called valid for \textit{polyhedron} $P\subseteq\Re^{n}$ if it is satisfied by all points of $P$.\end{definition}

\begin{definition}[\textbf{Face}]
If  $\pi x \leq \pi_{0}$ is a valid inequality for $P$ and $F=\{x\in P : \pi x=\pi_{0} \}$ then $F$ is called a face of $P$, 
represented by $(\pi,\pi_{0})$. ($F\neq \emptyset$ iff $max\{\pi_{x}:x \in P\}=\pi_{0}$).\end{definition}
\begin{definition}[\textbf{Facet}]\label{facet}
A  face of $P,$ is called a facet if $dim(F)=dim(P)-1.$
\end{definition}
\begin{prop}[\textbf{Description of P}]\label{necessary_suff}The facets of $P$ are necessary and sufficient for the description of $P$.
\end{prop}
\begin{definition}[\textbf{Convex Hull}]\label{convhull}
The set $S$ is defined as $S=P \cap Z^{n} $ where $P\subseteq\Re ^{n} $. The \textit{convex hull} of $S$ is denoted as $conv(S)$ and  $conv(S)=\{y \in \Re^{n}:y=\sum_{i=1}^{k}\lambda_{i}x^{i}\}$, where $x^{i} \in S$ and $\lambda_{i} \in \Re_{+} ~\forall i$, with $\sum_{i=1}^{k}\lambda_{i}=1.$

\end{definition}
\begin{figure}[hbp!]
\centering
\includegraphics[angle=270,
scale=0.65]{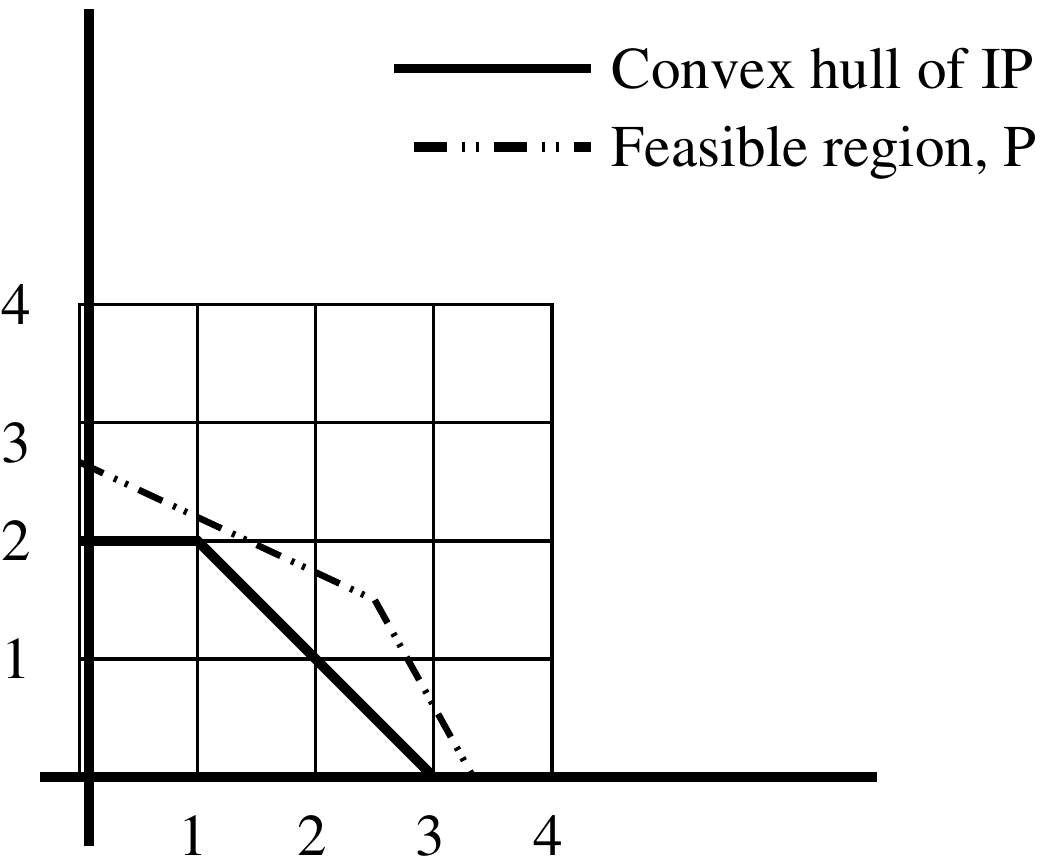}
\caption[Integral Polyhedron and its Convex Hull]{The dash-dotted line indicates the feasible set $P$, while the
solid line denotes the convex hull of the corresponding integral polyhedron on $\Re^{2}$.
An extremal point $y \in conv(S) $ has the property that $ \nexists \quad x^{1}, x^{2} \in conv(S)$ 
such that $y=\frac{1}{2}x^{1}+\frac{1}{2}x^{2}.$}
\label{extremal}
\end{figure}

\begin{prop}\label{dim_convS}
If $\pi_{x}\leq \pi_{0}$ defines a face of dimension $k-1$ of $conv(S)$, $\exists \quad k$ affinely independent points $x^{1}, ....,x^{k} \in S$ such that $\pi x^{i}=\pi_{0}$, for $i=1,...,k$.
\end{prop}
\begin{definition}[\textbf{Extreme point}]
A point $x\in P$ is called an extreme point if $\nexists x^{1}, x^{2} \in P$ such that $x=\frac{1}{2}x^{1}+\frac{1}{2}x^{2}.$
\end{definition}
\begin{prop}\label{facet_S_convS}
A maximal valid inequality for $S$, is one that dominates all other valid inequalities. The set of maximal valid inequalities of $S$ contains all the facet-defining inequalities of $conv(S)$, p.207 \cite{Wolsey}.
\end{prop}

The development in Wolsey~\cite{Wolsey}, uses a theorem relating the extremal points of the convex hull of $S$, to those
 of the set $S$ itself. They state a property without proof, for which the following intuitive proof can be given.
The point $y \in conv(S) $ is assumed to have the property that $\nexists x^{1}, x^{2} \in conv(S)$ 
such that $y=\frac{1}{2}x^{1}+\frac{1}{2}x^{2}$.
There are several cases depending on where the point $y$ lies in Figure~\ref{extremal}: 
\begin{enumerate}
\item $y \in conv(S)\backslash \{\textrm{boundary of }conv(S)\}$. This contradicts the assumption that $y$ is an extremal point of $conv(S)$.
This is because all strictly interior  points can be obtained as (non-trivial) convex combinations of some other points.
\item $y \in \{\textrm{boundary of~\textit{conv(S)}}\}\backslash \{x:x~\textrm{ is an extremal point of} S\}$. Refer to Figure~\ref{extremal}.
However, such points on the boundary of $conv(S)$ can be obtained as (non-trivial) convex combinations of the extremal points of $S$, (and $S \subseteq
 conv(S)$), thus contradicting the assumption. 
  \\
\end{enumerate}

 If $y$ does not lie in any of the regions mentioned above, then $y \in S$ and $S \subseteq conv(S)$. Now, if
$y$ is not an extremal point of $S$, $\exists \quad x^{1}, x^{2}$, such that $y=\frac{1}{2}x^{1}+\frac{1}{2}x^{2}$, where $x^{1},x^{2} \in
conv(S)$, contradicting the assumption again leading to the following lemma.

\begin{lemma}

 Hence, if such an extremal point of $conv(S)\quad \exists,$ it must be one of the extremal points of $S$.

\end{lemma}

\begin{theorem}\label{main}
If $P= \{x\in \Re^{n}: Ax \leq b\}$ is the feasible set of the LP relaxation of
the IP $max\{cx:x \in S\}$ where $ S=P\cap Z^{n}$ and
the LP on the convex hull of S, $max\{cx:x \in conv(S)\}$ has an optimal solution, $x^{\ast}$(which will be an extremal solution, since conv(S) is a polyhedron), then
$x^{\ast}$ will be an optimal solution to the IP on S too.

\end{theorem}
\subsection{Why the formulation used by Bertsimas et.al.,~\cite{Bert} is strong?}
We shall now try to explain the logic behind the occurrence of integral solutions when using the particularly strong formulation, given by Bertsimas et.al, \cite{Bert} using a constructive proof.

\subsubsection{Constructive Proof}

\begin{enumerate}

\item [(i)]Try to find the $conv(S)$, which is difficult usually. If we can find the conv(S) uniquely, we are done by~\textbf{Theorem~\ref{main}}. Else \ldots
\item[(ii)]Check which of the inequalities governing the IP are facets of the IP. This can be done by choosing each $a^{i}x\leq b_{i}$ in turn and finding the $dim(F)$ where $F$ is the set of points in IP satisfying $a^{i}x\leq b_{i}$ as  equality. If  $dim(F)=dim(IP)-1$, $F$ is a facet according to \textbf{Definition~\ref{facet}}.
 
\item[(iii)]The facets of $S$ are also facets of $conv(S)$ by~\textbf{Property~\ref{facet_S_convS}}. If, by using this iterative method, we can find the
 complete set
 of facets of $conv(S)$, we are done, as solving the LP on $conv(S)$ for an optimal solution (if it exists) will give an optimal solution to IP too, by~\textbf{Theorem \ref{main}}. However, even if we cannot get all the facets, of $conv(S)$ as is the case in real-life instances, the formulation is strong because of some constraints being facets and there is a high possibility of getting integral solutions.
 
\end{enumerate}
\section{Implementation }
We implemented the TFMP on artificially constructed flight data sets, to assess the computational performance of the model. Gnu Mathprog was used as the translator as it employs symbolic algebraic notation, which has advantages over languages like LINDO, which require that each constraint, alongwith the numerical values of the parameters be explicitly written down. The model was solved using the Open-source GLPK~\cite{Gnu} linear solver. The standalone solver contained within GLPK, \texttt{glpsol} was used for this project. 
\subsection{The Gnu Mathprog code}
Gnu Mathprog is an open-source code based on the C language. It is in fact the old version of AMPL as given in the paper by Kernighan et.al.~\cite{kerni}. The data consists of three primary sets:
\begin{enumerate}
\item The set of airports and sectors, denoted as $\mathscr{K}$. $$K=\{\textrm{Mumbai, Bangalore, Pune, Goa, Calicut,
Belgaum,Coimbatore}\}$$
\begin{figure}[hbp!]
\centering
\includegraphics[angle=0,
scale=0.6100]{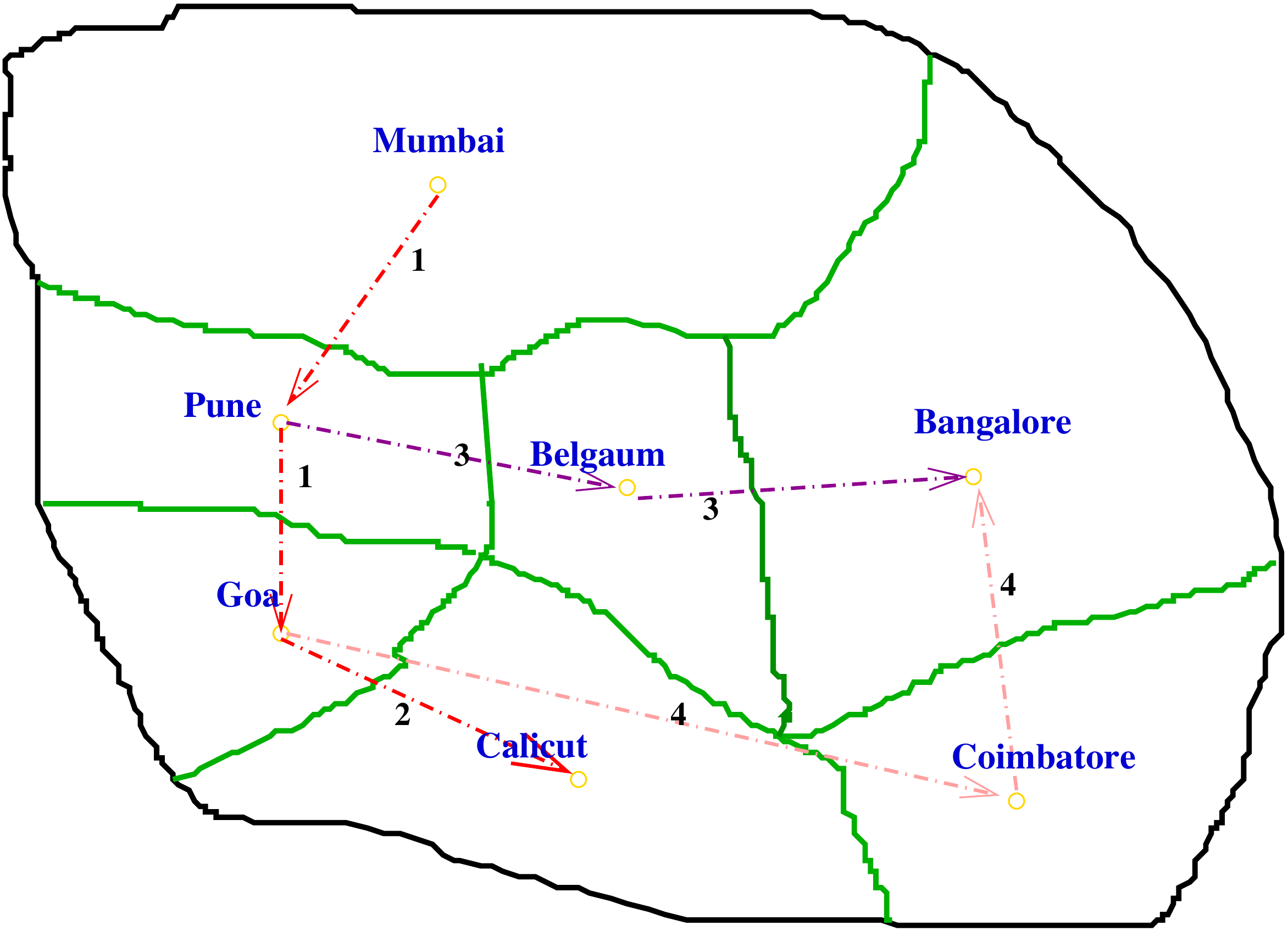}
\caption[Sectoral map]{The set of 7 flight sectors and 4 flights for South Western India. The Mumbai-Pune-Goa flight is continued by the Goa-Calicut flight. There are two other independent flights, Goa-Coimbatore-Bangalore and Pune-Belgaum-Bangalore. This is a hypothetical set of flights designed to illustrate concepts.}
\label{flights}
\end{figure}

\item The set of flights, denoted as $\mathscr{F}$.  $$K=\{\textrm{Mu\_Pu\_Go, Go\_Ca, Pu\_Be\_Ba, Go\_Co\_Ba
}\}$$ The Mu\_Pu\_Go flight is continued by the Go\_Ca flight. Other flights are independent of each other.
\item The set of time intervals, $\mathscr{T}$ over which the schedule is to be prepared.
\end{enumerate}
The Gnu Mathprog model for the TFMP is shown below. The model is mostly self-explanatory, but the reader is referred to Kernighan et.al ~\cite{kerni} for details.
\subsection{Gnu Mathprog Model}
\begin{verbatim}
/*
GMPL model for the TFMP.

Airports are:
Mumbai, Bangalore,Pune,Goa, Calicut, Belgaum, Coimbat
ore
Flights are:
* Mumbai-Pune-Goa connected to Goa-Calicut
*Pune-Belgaum-Bangalore
*Goa-Coimbatore-Bangalore
Time slots are t=0 to 10.

The chosen data set illustrates how the Pune-Belgaum-B
angalore flt.
is delayed because of restricted landing capacity at Ban
galore airport



*/




set N;
	#set of natural numbers
set K;
/*Airports or sectors*/
set F;
/*Flights*/
set T;
/*Time periods*/
 
param arr{k in K, t in T};
/*arrival capacity*/

param cap{k in K, t in T};
/* capacity*/

param atleast_one{f in F, k in K, t in T};
#flt. MUST reach sector k by time t

param dep{k in K, t in T};
/*departure capacity*/


param sched_dep{f in F};
/*scheduled dep*/

 param sched_arr{f in F};
/*sched arr time*/

/*Sector flight transit minimal time*/
param tr{f in F,k in K} ; 

param turn{f in F};
/*turnaround time for an aircraft*/

param ground{f in F};
/*ground holding cost per unit of time*/

param air{f in F};
/*air holding cost per unit of time*/

param feas{t in T, k in K, f in F};
		   #feasible flight-sector-time combinations
param connect{g in F, f in F};
				#set of connecting flights
param sec_ord {f in F , n in N} symbolic ;
	#represents the n th sector in flight f's  path		

param maxi{f in F};
		#cardinality of the set of sectors for a flight
#####VARIABLES######
var w{f in F, t in T , k in K}, >=0, <=1;
			#=1 if flt. f arrives at sector k BY time t
var grdel_fl{f in F}  ;
			#ground delay for flight f
var airdel_fl{f in F} ;
			#air delay for flight f

			#air delay for flight f
var cost;# cost variable =max(0, actual cost)
var dummyobj;#dummy objective function

minimize dummy:dummyobj;#the objective function

s.t. dummyobj1: dummyobj=sum{f in F}(ground[f]*grdel_fl[f]+air[f]*airdel_fl[f]);

s.t. pos:cost>=0;
s.t. pos1:cost>=dummyobj;
				# the objective function
s.t. grdel{f in F}:
grdel_fl[f]=sum{t in T:feas[t,sec_ord[f,1],f]=1
}t*(w[f,t,sec_ord[f,1]]-w[f,t-1,sec_ord[f,1]])-sched_dep[f];
						#ground delay for flight f

s.t. ardel{f in F}:
airdel_fl[f]=sum{t in T:feas[t,sec_ord[f, maxi[f]],f]=1 }
t*(w[f,t,sec_ord[f,maxi[f]]]-w[f,t-1,sec_ord[f, maxi[f]]])
-sched_arr[f]-grdel_fl[f];
						#air delay for flight f

s.t. domain{f in F,t in T, k in K:feas[t,k,f]=0.0}: 
w[f,t,k]=0;
		# no arrival/dep at infeasible combinations

s.t. atleast{f in F, k in K, t in T:atleast_one[f,k,t]=1}:
w[f,t,k]=1;
# cap on time at which flights must arrive at a sector

subject to depart {t in T, k in K } :
 sum{f in F:k=sec_ord[f,1] and feas[t,k,f]=1  }
(w[f,t,k]-w[f,t-1,k])<= dep[k,t];

#departure capacity of airport 
s.t. arrive {t in T,k in K} :
 sum{f in F:feas[t,k,f]=1 and k=sec_ord[f,maxi[f]]}
(w[f,t,k]-w[f,t-1,k])<=arr[k,t];
			#arrival capacity of airport
/*Sector  capacity*/
s.t. capa {t in T,k in K} : sum{f in F, n in N :n<maxi[f] and
feas[t,k,f]=1 and sec_ord[f,n]=k }
(w[f,t,sec_ord[f,n]]-w[f,t,sec_ord[f,n+1]])<=cap[k,t];

/*Sector connectivity*/
s.t. seccon {f in F,n in N, t in T: n<maxi[f] 
and feas[t,sec_ord[f,n],f]=1  } :
(w[f,t+tr[f,sec_ord[f,n]],sec_ord[f,n+1]]
-w[f,t,sec_ord[f,n]])<=0;

s.t. turnaround {f in F, g in F,t in T:
feas[t,sec_ord[f,1],f]=1 and
 connect[g,f]=1}:
(w[f,t,sec_ord[f,1]]-
w[g,t-turn[g],sec_ord[g,maxi[g]]])<=0;
				#turnaround constraint
/*defn of arrived by*/
s.t. arrvby {f in F,k in K,t in T: feas[t,k,f]=1 }:
(w[f,t,k]-w[f,t-1,k])>=0;
end;
\end{verbatim}
\subsection{Sample data for the model}
The data file can be maintained separately from the model, to permit different data to be solved using the same model.
The data has to be in a certain format for the interpreter-solver (\texttt{glpsol}) to work properly. Refer to \cite{Gnu} for details.
\subsubsection{The data file (unformatted)}
\begin{verbatim}
data;
# Mumbai Bangalore Pune Goa Calicut Belgaum Co
imbatore;
# MPG 1 GC 2 PbB 3 GCB 4
set K :=  Mu Ba Pu Go Ca Be Co;
set F :=Mu_Pu_Go Go_Ca Pu_Be_Ba Go_Co_Ba;
set T :=  0 1 2 3 4 5 6 7 8 9 10;
set N:= 1 2 3;

param maxi:=
Mu_Pu_Go 3
Go_Ca 2
Pu_Be_Ba 3
Go_Co_Ba 3;

param connect default 0.0 :=
[Mu_Pu_Go,Go_Ca]   1; 

param feas default 0.0 :=
[*,*,Mu_Pu_Go]:Mu  Ba     Pu    Go  Ca     Be     Co :=	  
0		        1	.	1	1	.	.	.	
1               	1      . 	1	1	.	.	.
2			1	.       1	1	.	.	.
3			1	.	1	1	.	.	.
4			.	.	1	1	.	.	.
5			.	.	.	1	.	.	.
6			.	.	.	.	.	.	.
7			.	.	.	.	.	.	.
8			.	.	.	.	.	.	.
9			.	.	.	.	.	.	.
10			.	.	.	.	.	.	.

[*,*,Pu_Be_Ba]:Mu	Ba	Pu	Go	Ca	Be	Co :=
0			.	.	.	.	.	.	.
1			.	.	.	.	.	.	.
2			.	.	1	.	.	.	.
3			.	.	1	.	.	1	.
4			.	.	1	.	.	1	.
5			.	1	1	.	.	1	.
6			.	1	.	.	.	1	.
7			.	1	.	.	.	1	.
8			.	.	.	.	.	.	.
9			.	.	.	.	.	.	.
10			.	.	.	.	.	.	.

[*,*,Go_Ca]:	Mu	Ba	Pu	Go	Ca	Be	Co :=
0			.	.	.	.	.	.	.
1			.	.	.	.	.	.	.
2			.	.       .	.	.	.	.
3			.	.	.	1	.	.	.
4			.	.	.	1	1	.	.
5			.	.	.	1	1	.	.	
6			.	.	.	1	1	.	.
7			.	.	.	.	1	.	.
8			.	.	.	.	1	.	.
9			.	.	.	.	.	.	.
10			.	.	.	.	.	.	.

[*,*,Go_Co_Ba]:Mu	Ba	Pu	Go	Ca	Be	Co :=
0			.	.	.	.	.	.	.
1			.	.	.	.	.	.	.
2			.	.       .	.	.	.	.
3			.	.	.	1	.	.	.
4			.	.	.	1	.	.	1
5			.	1	.	1	.	.	1	
6			.	1	.	.	.	.	1	
7			.	1	.	.	.	.	.
8			.	1	.	.	.	.	.
9			.	.	.	.	.	.	.
10			.	.	.	.	.	.	.;


param sec_ord :=
[Mu_Pu_Go,1]	Mu
[Mu_Pu_Go,2]	Pu
[Mu_Pu_Go,3]	Go
[Go_Ca,1]	Go
[Go_Ca,2]	Ca
[Pu_Be_Ba,1]	Pu
[Pu_Be_Ba,2]	Be
[Pu_Be_Ba,3]	Ba
[Go_Co_Ba,1]	Go
[Go_Co_Ba,2]	Co
[Go_Co_Ba,3]	Ba;





param arr default 2.0 :=	 
;

param dep default 2.0 := 
;

param cap default 2.0 :=
;

param tr:	               Mu	Ba	Pu	Go	Ca	Be	Co :=


Mu_Pu_Go		1	1	1	1	1	1	1
Go_Ca			1	1	1	1	1	1	1
Pu_Be_Ba		        1	1	1	1	1	1	1
Go_Co_Ba		1	1	1	1	1	1      1;

param turn[Mu_Pu_Go]:=1;

param sched_dep:= 

Mu_Pu_Go 1
Go_Ca    4
Pu_Be_Ba 3 
Go_Co_Ba 3;

param sched_arr:= 
Mu_Pu_Go 3
Go_Ca 	 5
Pu_Be_Ba 5
Go_Co_Ba 5;

param air := 
Mu_Pu_Go 20000
Go_Ca 	 30000
Pu_Be_Ba 40000
Go_Co_Ba 50000;

param ground:= 
Mu_Pu_Go 400
Go_Ca 	 600
Pu_Be_Ba 800
Go_Co_Ba 1000;

param atleast_one default 0.0 :=

[Mu_Pu_Go,*,*]	Mu	3	1
		                Pu	4	1
		                Go	5	1
[Go_Ca,*,*]		Go	6	1
				Ca	7	1
[Pu_Be_Ba,*,*]		Pu	4	1
				Be	6	1
				Ba	7	1
[Go_Co_Ba,*,*]		Go	4	1
				Co	6	1
				Ba	7	1;		

end;
\end{verbatim}
\section{Solution}
The LP relaxation of the TFMP model is first solved with  $\forall w_{f,t}^{k} \in \Re \cap [0,1]$ to verify if our data set supports the claims made in the paper by Bertsimas et.al ~\cite{Bert}.
\subsection{Solution to the LP relaxation of TFMP with conflicting flights}
A set of  imaginary  flights Pu\_Be\_Ba and Go\_Co\_Ba having the same scheduled arrival timings at Bangalore, which has an arrival capacity for only one flight  was created, as in Table~\ref{tab_conflict}. However, even for this conflicting data set, the LP relaxation still had integral solutions, lending support to the results due to Bertsimas et.al. The Pu\_Be\_Ba flight is held for one interval after its scheduled departure from Goa, in order to satisfy the conflicting constraints. The minimal cost is 800 corresponding to the value of the Pu\_Be\_Ba flight, $c^{g}_{Pu\_Be\_Ba}$.
\begin{table}[hbp!]
\begin{tabular}{ccccccc}
\toprule
\multirow{2}{*}{\textbf{Flight}}&\multicolumn{2}{c}{\textbf{Departure time}} &\multicolumn{2}{c}{\textbf{Arrival time}}&\multirow{2}{*}{\textbf{Ground delay}}
&\multirow{2}{*}{\textbf{Air delay}}\\
 \cmidrule{2-5}
&\textbf{Scheduled}&\textbf{Actual}&\textbf{Scheduled}&\textbf{Actual}&&\\\midrule
Mu\_Pu\_Go & 1 & 1& 3 &3&0&0\\
Go\_Ca & 4 &4&5&5&0&0\\
Pu\_Be\_Ba &3 &4&$5^{\star} $&6&1&0\\
Go\_Co\_Ba &3 &3&$5^{\star} $&5&0&0\\

\midrule
\multicolumn{7}{c}{All airports have arrival capacity =1, at all  times}\\

\bottomrule
\end{tabular}
\caption[Scenario 1: Conflicting flights]{Solution to the LP relaxation of TFMP with conflicting flights.\\ 
\small{ \textit{Conflicting data set with Pu\_Be\_Ba and Go\_Co\_Ba having the same scheduled arrival timings (as shown by the asteriks ($\star$)) at Bangalore, which can accomodate only one arrival at all time intervals. Notice in the following solution, how the  Pu\_Be\_Ba flight suffers a ground delay of 1 unit to resolve the conflict.}}}
\label{tab_conflict}

 \end{table}

\begin{verbatim}
Problem:    lp_tfmp
Rows:       588
Columns:    316
Non-zeros:  582
Status:     OPTIMAL
Objective:  cost = 800 (MINimum)

   No.   Row name   St   Activity     Lower bound   Upper bound    Marginal
------ ------------ -- ------------- ------------- ------------- -------------
     1 cost         B            800                             
             

   No. Column name  St   Activity     Lower bound   Upper bound    Marginal
------ ------------ -- ------------- ------------- ------------- -------------
     1 w[Mu_Pu_Go,0,Mu]
                    B              0             0             1 
     2 w[Mu_Pu_Go,1,Mu]
                    B              1             0             1 
     3 w[Mu_Pu_Go,2,Mu]
                    B              1             0             1 
     4 w[Mu_Pu_Go,3,Mu]
                    B              1             0             1 
     5 w[Go_Ca,2,Go]
                    B              0             0             1 
     6 w[Go_Ca,3,Go]
                    B              0             0             1 
     7 w[Go_Ca,4,Go]
                    B              1             0             1 
     8 w[Go_Ca,5,Go]
    
   309 grdel_fl[Mu_Pu_Go]
                    B              0             0               
   310 grdel_fl[Go_Ca]
                    B              0             0               
   311 grdel_fl[Pu_Be_Ba]
                    B              1             0               
   312 grdel_fl[Go_Co_Ba]
                    B              0             0               
   313 airdel_fl[Mu_Pu_Go]
                    NL             0             0                       19600
   314 airdel_fl[Go_Ca]
                    B              0             0               
   315 airdel_fl[Pu_Be_Ba]
                    B              0             0               
   316 airdel_fl[Go_Co_Ba]
                    B              0             0               

Karush-Kuhn-Tucker optimality conditions:

KKT.PE: max.abs.err. = 0.00e+00 on row 0
        max.rel.err. = 0.00e+00 on row 0
        High quality

KKT.PB: max.abs.err. = 0.00e+00 on row 0
        max.rel.err. = 0.00e+00 on row 0
        High quality

KKT.DE: max.abs.err. = 0.00e+00 on column 0
        max.rel.err. = 0.00e+00 on column 0
        High quality

KKT.DB: max.abs.err. = 0.00e+00 on row 0
        max.rel.err. = 0.00e+00 on row 0
        High quality

End of output
\end{verbatim}
\subsection{Baseline solution without conflicts} 
The LP relaxation of the TFMP is now solved with the baseline data set, as in Table~\ref{tab_base}.
As can be seen from the solution below, the linear program adjusts the departure and arrival timings so that the total cost is minimum. For the baseline  data set, the cost is zero, implying that there are no delays of any sort.
\begin{table}[hbp!]
\begin{tabular}{ccccccc}\toprule
\multirow{2}{*}{\textbf{Flight}}&\multicolumn{2}{c}{\textbf{Departure time}} &\multicolumn{2}{c}{\textbf{Arrival time}}&\multirow{2}{*}{\textbf{Ground delay}}
&\multirow{2}{*}{\textbf{Air delay}}\\
 \cmidrule{2-5}
&\textbf{Scheduled}&\textbf{Actual}&\textbf{Scheduled}&\textbf{Actual}&&\\\midrule
Mu\_Pu\_Go & 1 & 1& 3 &3&0&0\\
Go\_Ca & 4 &4&5&5&0&0\\
Pu\_Be\_Ba &3 &3&$5^{\star}$&5&0&0\\
Go\_Co\_Ba &3 &3 &$5^{\star}$&5&0&0\\

\midrule
\multicolumn{7}{c}{All airports have arrival capacity =2, at all  times}\\

\bottomrule
\end{tabular}
\caption[Scenario 2: Baseline solution]{Baseline solution without conflicts.\\ \small{ \textit{Note the baseline data set is conflict free, since though Pu\_Be\_Ba and Go\_Co\_Ba have the same scheduled arrival timings (as shown by the asteriks ($\star$)) at Bangalore, both can be accomodated due to the increased arrival capacity of 2. Notice that all filghts are delay-free, in this instance.}}}
\label{tab_base}

 \end{table}

\begin{verbatim}
Problem:    tfmp
Rows:       588
Columns:    316
Non-zeros:  582
Status:     OPTIMAL
Objective:  cost = 0 (MINimum)

   No.   Row name   St   Activity     Lower bound   Upper bound    Marginal
------ ------------ -- ------------- ------------- ------------- -------------
     1 cost         B              0                             
    

   No. Column name  St   Activity     Lower bound   Upper bound    Marginal
------ ------------ -- ------------- ------------- ------------- -------------
     1 w[Mu_Pu_Go,0,Mu]
                    B              0             0             1 
     2 w[Mu_Pu_Go,1,Mu]
                    B              1             0             1 
     3 w[Mu_Pu_Go,2,Mu]
                    B              1             0             1 
     4 w[Mu_Pu_Go,3,Mu]
                    B              1             0             1 
     5 w[Go_Ca,2,Go]
                    B              0             0             1 
     6 w[Go_Ca,3,Go]
                    B              0             0             1 
     7 w[Go_Ca,4,Go]
                    B              1             0             1 
     8 w[Go_Ca,5,Go]
                    NU             1             0             1          -600
    
   309 grdel_fl[Mu_Pu_Go]
                    B              0             0               
   310 grdel_fl[Go_Ca]
                    B              0             0               
   311 grdel_fl[Pu_Be_Ba]
                    B              0             0               
   312 grdel_fl[Go_Co_Ba]
                    B              0             0               
   313 airdel_fl[Mu_Pu_Go]
                    NL             0             0                       19600
   314 airdel_fl[Go_Ca]
                    B              0             0               
   315 airdel_fl[Pu_Be_Ba]
                    B              0             0               
   316 airdel_fl[Go_Co_Ba]
                    B              0             0               

Karush-Kuhn-Tucker optimality conditions:

KKT.PE: max.abs.err. = 0.00e+00 on row 0
        max.rel.err. = 0.00e+00 on row 0
        High quality

KKT.PB: max.abs.err. = 0.00e+00 on row 0
        max.rel.err. = 0.00e+00 on row 0
        High quality

KKT.DE: max.abs.err. = 0.00e+00 on column 0
        max.rel.err. = 0.00e+00 on column 0
        High quality

KKT.DB: max.abs.err. = 0.00e+00 on row 0
        max.rel.err. = 0.00e+00 on row 0
        High quality

End of output
\end{verbatim}
\subsection{Arrivals/Departures before scheduled time: negative costs}
A further interesting instance is  when the actual arrival/departure timings are before the scheduled timings.
Clearly, there is nothing restrictive about the formulation which will prevent such an instance from being considered. However, the economic interpretation is not the same. Moreover, with the objective being of minimization, the solver forces high values of negative air delays. In fact, there might be a situation, where one flight has a very high negative air delay and others have low positive values of ground delay. The real total cost, then is not negative. Actually, wherever a cost variable has a negative value it can be substituted by an auxiliary variable indicating its real value as $$\alpha \geq \max(0,\beta),$$
where Real cost$=\alpha$ and computed cost=$\beta$. The objective function will be $\beta$ while the problem is to minimize it.
The  solution sumarized in Table~\ref{tab_negcost}, indicates such an instance of a negative cost, due to a negative air delay  for the Go\_Ca\_\_Ba flight.
\begin{table}[hbp!]
\begin{tabular}{ccccccc}\toprule
\multirow{2}{*}{\textbf{Flight}}&\multicolumn{2}{c}{\textbf{Departure time}} &\multicolumn{2}{c}{\textbf{Arrival time}}&\multirow{2}{*}{\textbf{Ground delay}}
&\multirow{2}{*}{\textbf{Air delay}}\\
 \cmidrule{2-5}
&\textbf{Scheduled}&\textbf{Actual}&\textbf{Scheduled}&\textbf{Actual}&&\\\midrule
Mu\_Pu\_Go & 1 & 1& 3 &3&0&0\\
Go\_Ca & 4 &4&5&5&0&0\\
Pu\_Be\_Ba &3 &3&5&5&0&0\\
Go\_Co\_Ba &3 &$1^{\star} $&5&3&-2&0\\

\midrule
\multicolumn{7}{c}{All airports have arrival capacity =2, at all  times}\\

\bottomrule
\end{tabular}
\caption[Scenario 3: Arrivals/Departures before-time]{Arrivals/Departures before scheduled time: negative costs.\\ \small{ \textit{The feasible set of times $T_{\textrm{Go\_Co\_Ba}}$ is extended backward in time by 2 intervals, for all airports and sectors for the
Go\_Co\_Ba flight. This permits the flight to depart before time from Goa airport, (see asterik ($\star$)) and enter all sectors before its scheduled 
times, finally arriving before-time at Bangalore airport. In particular, note how the entire freedom-to-adjust (of 2 intervals) available in the schedule is taken
up by the  solution.}}}
\label{tab_negcost}

 \end{table}

\begin{verbatim}
Problem:    tfmp
Rows:       591
Columns:    318
Non-zeros:  594
Status:     OPTIMAL
Objective:  dummy = -2000 (MINimum)


   No. Column name  St   Activity     Lower bound   Upper bound    Marginal
------ ------------ -- ------------- ------------- ------------- -------------
     1 w[Mu_Pu_Go,0,Mu]
                    B              0             0             1 
     2 w[Mu_Pu_Go,1,Mu]
                    B              1             0             1 
     3 w[Mu_Pu_Go,2,Mu]
                    NU             1             0             1          -400
     4 w[Mu_Pu_Go,3,Mu]
                    NU             1             0             1        -58800
     5 w[Go_Ca,3,Go]
                    B              0             0             1 
      309 grdel_fl[Mu_Pu_Go]
                    B              0                             
   310 grdel_fl[Go_Ca]
                    B              0                             
   311 grdel_fl[Pu_Be_Ba]
                    B              0                             
   312 grdel_fl[Go_Co_Ba]
                    B             -2                             
   313 airdel_fl[Mu_Pu_Go]
                    B              0                             
   314 airdel_fl[Go_Ca]
                    B              0                             
   315 airdel_fl[Pu_Be_Ba]
                    B              0                             
   316 airdel_fl[Go_Co_Ba]
                    B              0                             
   317 cost         B              0                             
   318 dummyobj     B          -2000                             

Karush-Kuhn-Tucker optimality conditions:

KKT.PE: max.abs.err. = 5.68e-12 on row 4
        max.rel.err. = 2.27e-13 on row 2
        High quality

KKT.PB: max.abs.err. = 4.44e-16 on column 14
        max.rel.err. = 2.22e-16 on row 579
        High quality

KKT.DE: max.abs.err. = 6.26e-10 on column 16
        max.rel.err. = 3.71e-10 on column 305
        High quality

KKT.DB: max.abs.err. = 0.00e+00 on row 0
        max.rel.err. = 0.00e+00 on row 0
        High quality

End of output
\end{verbatim}
\subsection{Departure of an outgoing flight scheduled before arrival of the incoming connecting flight}
This solution, seen in Table~\ref{tab_connect}, indicates how the Go\_Ca flight is delayed on the ground at Goa airport, because its scheduled departure is before the
scheduled arrival of its incoming-connecting Mu\_Pu\_Go flight at Goa airport.
\begin{table}[hbp!]
\begin{tabular}{ccccccc}\toprule
\multirow{2}{*}{\textbf{Flight}}&\multicolumn{2}{c}{\textbf{Departure time}} &\multicolumn{2}{c}{\textbf{Arrival time}}&\multirow{2}{*}{\textbf{Ground delay}}
&\multirow{2}{*}{\textbf{Air delay}}\\
 \cmidrule{2-5}
&\textbf{Scheduled}&\textbf{Actual}&\textbf{Scheduled}&\textbf{Actual}&&\\\midrule
Mu\_Pu\_Go & 1 & 1& $3^{\star}$ &3&0&0\\
Go\_Ca &$2^{\star}$ &4&4&6&2&0\\
Pu\_Be\_Ba &3 &3&5&5&0&0\\
Go\_Co\_Ba &3 &3&5&5&0&0\\

\midrule
\multicolumn{7}{c}{All airports have arrival capacity =2, at all  times}\\

\bottomrule
\end{tabular}
\caption[Scenario 4: Inconsistent timings]{Departure of an outgoing flight scheduled before arrival of the incoming connecting flight.\\\small{ \textit{Here, at Goa airport, the departing flight Go\_Ca, is scheduled to depart, before the scheduled arrival of its incoming, connecting flight Mu\_Pu\_Go, arriving from Mumbai. The asteriks ($\star$) indicate a situation where passengers getting off  the Mu\_Pu\_Go flight at Goa airport would miss their connecting flight, Go\_Ca, to Calicut. Notice how the LP resolves this issue by holding the  Go\_Ca flight on the ground for 2 time intervals, beyond its scheduled departure, so as to accomodate passengers from the incoming flight. Note that the schedule and connectivity parameters are such, that this delay is inevitable, and the LP ensures that exactly those flights are delayed which minimize the cost of delay. }}}
\label{tab_connect}

 \end{table}

\begin{verbatim}
Problem:    tf
Rows:       585
Columns:    316
Non-zeros:  560
Status:     OPTIMAL
Objective:  cost = 1200 (MINimum)

   
   No. Column name  St   Activity     Lower bound   Upper bound    Marginal
------ ------------ -- ------------- ------------- ------------- -------------
     1 w[Mu_Pu_Go,0,Mu]
                    B              0             0             1 
     2 w[Mu_Pu_Go,1,Mu]
                    B              1             0             1 
     3 w[Mu_Pu_Go,2,Mu]
                    B              1             0             1 
     4 w[Mu_Pu_Go,3,Mu]
                    B              1             0             1 
     5 w[Go_Ca,3,Go]
                    B              0             0             1 
     6 w[Go_Ca,4,Go]
                    B            0.5             0             1 
     7 w[Go_Ca,5,Go]
                    B            0.5             0             1 
     8 w[Go_Ca,6,Go]
                    B              1             0             1 
     9 w[Pu_Be_Ba,2,Pu]
                    B              0             0             1 
    10 w[Pu_Be_Ba,3,Pu]
                    B              1             0             1 
   309 grdel_fl[Mu_Pu_Go]
                    NL             0             0                        1000
   310 grdel_fl[Go_Ca]
                    B              2             0               
   311 grdel_fl[Pu_Be_Ba]
                    B              0             0               
   312 grdel_fl[Go_Co_Ba]
                    B              0             0               
   313 airdel_fl[Mu_Pu_Go]
                    NL             0             0                       20600
   314 airdel_fl[Go_Ca]
                    NL             0             0                       30000
   315 airdel_fl[Pu_Be_Ba]
                    B              0             0               
   316 airdel_fl[Go_Co_Ba]
                    B              0             0               

Karush-Kuhn-Tucker optimality conditions:

KKT.PE: max.abs.err. = 0.00e+00 on row 0
        max.rel.err. = 0.00e+00 on row 0
        High quality

KKT.PB: max.abs.err. = 0.00e+00 on row 0
        max.rel.err. = 0.00e+00 on row 0
        High quality

KKT.DE: max.abs.err. = 0.00e+00 on column 0
        max.rel.err. = 0.00e+00 on column 0
        High quality

KKT.DB: max.abs.err. = 0.00e+00 on row 0
        max.rel.err. = 0.00e+00 on row 0
        High quality

End of output
\end{verbatim}

\section{Conclusions }We have tried to explore the logic behind the claims by Bertsimas et.al about integral solutions to the LP relaxation of the TFMP.

Polyhedral theory only indicates that the stronger TFMP formulation of Bertsimas et.al might lead to integral solutions in some cases.
Unless we obtain all facet-defining inequalities of the
\textit{conv(S)}, we are not assured of obtaining integral
solutions from the LP relaxation. Therefore, only computations can provide a rough estimate of the frequency of integral solutions, in the worst case, which in turn will decide whether the formulation offers practical advantages. Our computations indicate that the encouraging results reported by Bertsimas et.al are not merely fortuitous or due
to their specific data set. Indeed, we found that the TFMP had integral solutions   even in case of artificial data sets generated to include severe conflicts in the flight schedules. In our limited tests with 4-5 scenarios, we obtained non-integral solutions only once.

This is of significant practical importance because, the LP relaxation can be solved even on small machines with
  low memory and processor speed. The formulation  with 308 variables took less than 0.1 seconds when solved on a regular IBM laptop having 1GB RAM and 1.66 GHz Processor speed. In contrast, the Integer program has to be solved by Branch and Bound which would be much more expensive and difficult to  implement on such a machine, especially since the number of variables rises rapidly with addition of new flights. The importance of obtaining integral solutions, using the LP relaxation itself (without the need for Branch and Bound), will be more acutely felt for problems of large sizes. 
\section{Directions for future work}
An important observation, is that the size of the problem is quite large even for simple cases. E.g. for a data set with 4 flights, 7 sectors and 11 time slots, the formulation already has 308 variables. In a real life data set, the proportion of conflicting flights may not be very high. In such a situation, it is better to solve the LP for only those flights which are in conflict with 
each other. 

One approach could be, to use the scheduled arrival and departure timings,to generate a set of flights which are in conflict with each other, called the set$ \mathscr{X}$. The set $\mathscr{Y}=\mathscr{F}\backslash \mathscr{X}$ then consists of flights which are not in conflict either among themselves or with those in $\mathscr{X}$ . The TFMP is solved only for  flights $f\in \mathscr{X}$. Then, the entire set $\mathscr{F}$ is again checked for conflicts, with the flight timings for those in $\mathscr{X}$ now updated to those given by the TFMP solution.
If this causes some flights in  $\mathscr{Y}$ to be in
conflict with those in  $\mathscr{X}$, such flights
are transferred from $\mathscr{Y}$ to $\mathscr{X}$
. This iterative process is continued till we get a schedule with no conflicts.\\The problem with such an approach is that
the set $\mathscr{X}$ might grow, until finally $\mathscr{X} \equiv \mathscr{F}$ thus offering no advantages in terms of  problem  size. Moreover, the TFMP has to be solved once per iteration and further analysis needs to be done to weigh the benefits of solving a reduced problem, against the expense in solving the TFMP several times. The concept, is however very attractive for an initial schedule having very few conflicting flights, if some heuristics can be developed to locally adjust the conflicting set $\mathscr{X}$ without affecting the larger, non-conflicting set $\mathscr{Y}$.

\end{document}